\documentclass[11pt,letterpaper,twoside]{amsart}
\usepackage[top=.8in, bottom=.8in, left=.8in, right=.8in]{geometry}
\usepackage[utf8]{inputenc}
\usepackage{amsmath}
\usepackage{amssymb}
\usepackage{amsfonts}
\usepackage{amsthm}
\usepackage{graphicx}
\usepackage{color}
\usepackage[font=footnotesize,labelfont=bf]{caption}
\usepackage{hyperref}
\usepackage[nameinlink]{cleveref}
\usepackage[alphabetic]{amsrefs}
\usepackage{titlesec}
\usepackage{titlecaps}

\theoremstyle{definition}
\newtheorem{theorem}{Theorem}

\newtheorem{defn}[theorem]{Definition}
\newtheorem{remark}[theorem]{Remark}

\title{\textbf{SPLITTING SPHERES FOR UNLINKED $S^2$'S IN $S^4$}} 
\author{ALISON TATSUOKA}
\thanks {The author was partially supported by a NSF Graduate Research Fellowship.}
\address{Department of Mathematics, Princeton University, Princeton, NJ 08540, USA}
\email{at8451@princeton.edu}
\date{}

\begin{document}

\maketitle

\noindent \textbf{Abstract.} We show that there exist infinitely many pairwise non-isotopic splitting spheres for two unlinked, unknotted $S^2$'s in $S^4$. This answers a question posed by Hughes, Kim, and Miller.

\titleformat{\section}[block]{\filcenter\scshape}{\thesection.}{0.5em}{\titlecap}
 \section{Introduction}

\noindent In this paper, we study \textit{splitting spheres} for the unlink of two unknotted $S^2$'s in $S^4$.

\begin{defn} \label{definition of splitting sphere}
    Let $K=K_1\sqcup K_2$ be a 2-component link in $S^n$. A \textit{splitting sphere} for $K$ is an embedded $S^{n-1}\subset S^n \setminus K$ such that $K_1$ and $K_2$ lie in distinct components of $S^n \setminus K$. Two splitting spheres are isotopic if they are isotopic in $S^n \setminus K$.
\end{defn}

\noindent In dimension 3, any two splitting spheres for a 2-component link are isotopic. In dimension 4, however, Hughes, Kim, and Miller \cite{HKM} showed that if $K$ is the unlink of a genus $g$ and a genus $h$ surface, where $g\geq 4$, then there exists a pair of non-isotopic splitting spheres for $K$; they also showed the same result for the cases $(g,h)=(2,2), (2,3)$, and $(3,3)$. Noting that their methods didn't extend to lower genus unlinks, they asked whether there exist multiple splitting spheres when $(g,h)=(0,0)$, that is, when $K$ is the unlink of two unknotted $S^2$'s in $S^4$. We give an affirmative answer to their question by providing an infinite family of pairwise non-isotopic $S^3$'s in $S^4 \setminus K \approx S^1 \times B^3 \# S^1 \times B^3$.

\begin{theorem} \label{splitting}
    Let $K$ be the unlink of two unknotted $S^2$'s in $S^4$. Then there exist infinitely many pairwise non-isotopic splitting spheres for $K$.
\end{theorem}

\begin{remark}
    One reason splitting spheres might be of interest is in light of the 4-dimensional smooth Schoenflies Conjecture, which asks if there exist non-isotopic $S^3$'s in $S^4$. Budney-Gabai \cite{BG} and independently Watanabe \cite{Wa} showed that there exist infinitely many $B^3$'s in $S^1 \times B^3$ with boundary $pt \times \partial B^3$ which are pairwise non-isotopic rel boundary; as a consequence, the examples of \cite{BG} provide infinitely many $B^3$'s in $S^4$ which are non-isotopic rel boundary. In the same paper \cite{BG} also construct infinitely many pairwise non-isotopic $S^3$'s in $S^1 \times S^3$. Other examples of knotted $S^3$'s come from Iida-Konno-Mukherjee-Taniguchi \cite{IKMT} and Konno-Mukherjee-Taniguchi \cite{KMT}, who found infinitely many non-isotopic $S^3$'s in 4-manifolds with nontrivial second homology (the simplest being $S^2 \times S^2 \# S^2 \times S^2$), and from the splitting spheres in surface link complements of Hughes-Kim-Miller in \cite{HKM}. We note that none of their methods extend to the case of Theorem \ref{splitting}.
\end{remark}

\begin{figure}[h!]
    \centering
    \includegraphics[width=0.3\linewidth]{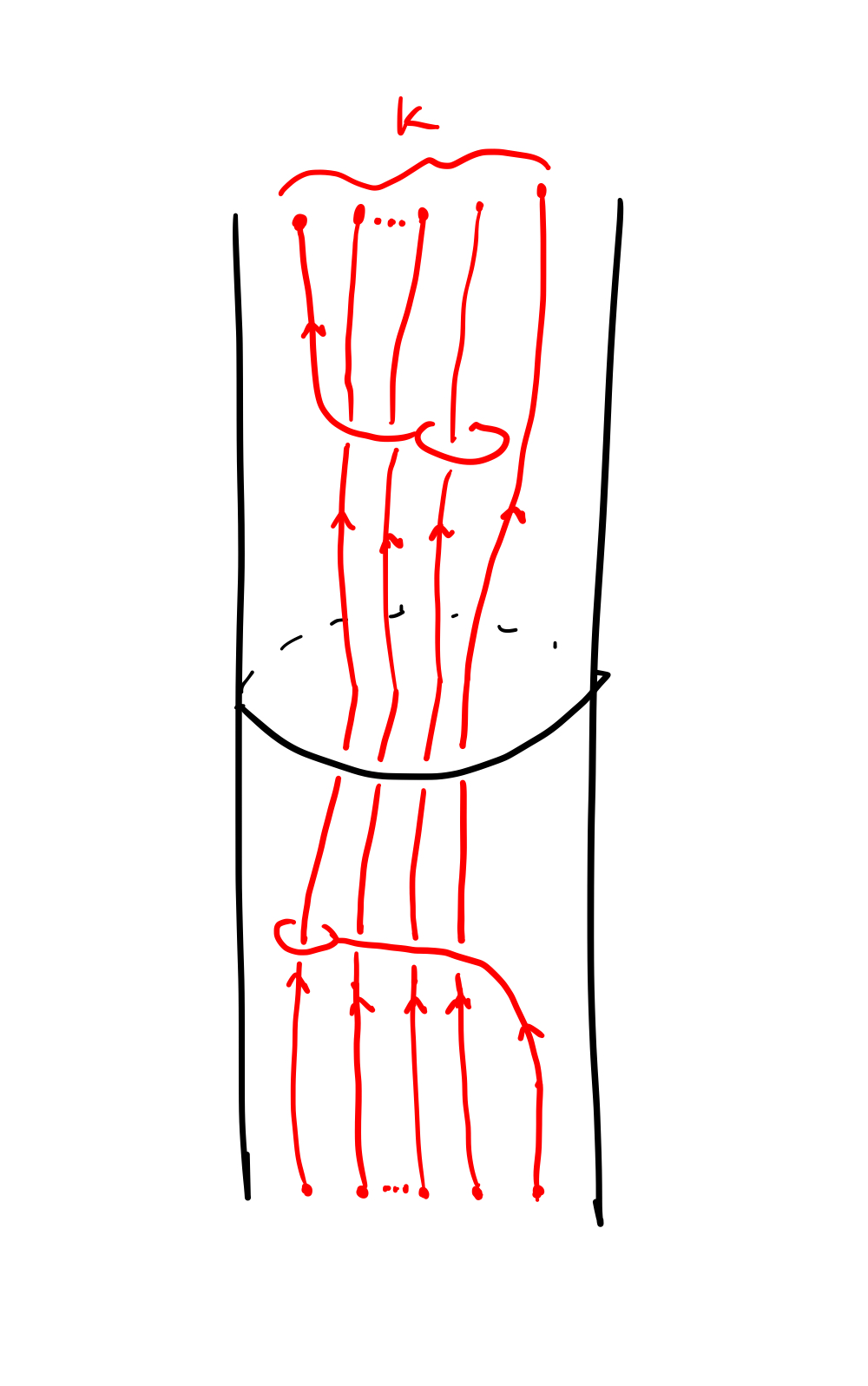}
    \caption{The barbell representing the diffeomorphism $\delta_k$ of $S^1 \times B^3$, as seen in the $S^1 \times D^2 \times 0$ slice \cite{BG}.}
    \label{delta}
\end{figure}

\noindent Our result uses the work of Budney-Gabai in \cite{BG}, where they find knotted 3-balls in $S^1 \times B^3$ by constructing an infinite family of barbell diffeomorphisms $\delta_k$ of $S^1 \times B^3$ which are pairwise non-isotopic rel boundary and nontrivial for $k \geq 4$ in $\pi_0(\text{Diff}_\partial(S^1 \times B^3)/\text{Diff}_\partial(B^4))$. We quickly recall the definition of a barbell diffeomorphism here; for more details see Chapter 5 of \cite{BG}. A \textit{barbell diffeomorphism} of a 4-manifold $X$ is a diffeomorphism of $X$ defined by embedding an $S^2 \times D^2 \natural S^2 \times D^2$ into the interior of $X$ and pushing forward the following diffeomorphism of $S^2 \times D^2 \natural S^2 \times D^2$: let $A_1$ and $A_2$ be two disjoint, properly embedded arcs in $B^4$. Spin $A_1$ about $A_2$ keeping $\partial A_1$ fixed; this isotopy extends to a diffeomorphism of $B^4$ that fixes neighborhoods of $A_1$ and $A_2$, and deleting those neighborhoods provides the desired diffeomorphism of $S^2 \times D^2 \natural S^2 \times D^2 \approx B^4 \setminus A_1 \sqcup A_2$. The embedding $S^2 \times D^2 \natural S^2 \times D^2 \hookrightarrow X$ is called the \textit{barbell implantation} and is encoded pictorially as a \textit{barbell}; the barbell diffeomorphism is supported in a neighborhood of the barbell. Viewing $S^1 \times B^3$ as $S^1\times D^2 \times [-1,1]$, \Cref{delta} depicts the barbell in $S^1 \times B^3$ corresponding to the $\delta_k$ diffeomorphism of \cite{BG}, as seen in the $S^1 \times D^2 \times 0$ slice. 

\noindent \textbf{Conventions.} We work in the smooth category. We denote by $\text{Diff}_\partial (X)$ the group of diffeomorphisms of $X$ which are the identity on $\partial X$. We consider the quotient $\text{Diff}_\partial (X) / \text{Diff}_\partial (B^4)$, where two diffeomorphisms are equivalent if they agree outside of a $B^4 \subset X$, and the group $\pi_0(\text{Diff}_\partial (X) / \text{Diff}_\partial (B^4))$, where two diffeomorphisms are equivalent if they are isotopic rel $\partial$ to diffeomorphisms that agree outside of a $B^4 \subset X$. We use $\approx$ to denote isotopy and all isotopies are rel boundary where applicable. We sometimes write $f\approx g$ rel $B^4$ to mean that $f$ and $g$ are equivalent in $\pi_0(\text{Diff}_\partial (X) / \text{Diff}_\partial (B^4))$.

\noindent \textbf{Acknowledgements.} The author would like to thank Seungwon Kim for introducing her to this problem, her advisor Dave Gabai for his comments and encouragement, and Gheehyun Nahm for helpful conversations.

\section{Proof of Theorem 2}
\begin{proof}

    Let $L \# R$ denote $S^1 \times B^3 \# S^1 \times B^3$, where $L$ and $R$ are $S^1 \times B^3$'s. Let $\Sigma$ denote the $S^3$ along which $L$ and $R$ are identified in the connect sum.

    \noindent Consider the barbell diffeomorphism $\beta_k: L \# R \rightarrow L \# R$ given by the red barbell implantation below, $k\geq 4$. As in \cite{BG}, we think of $S^1 \times B^3$ as $S^1 \times D^2 \times [-1,1]$, where the $t\in[-1,1]$ parameter represents time. We think of $L \# R$ as two copies of $(S^1 \times D^2) \times [-1,1] \setminus D^3 \times (-0.5, 0.5)$ ($L$ on the left and $R$ on the right) which are identified along the $S^3$ boundary of the $B^4=D^3 \times [-0.5,0.5]$. \Cref{beta} shows the barbell in the $t=0$ slice of $L\# R$. The orange $S^2$'s on left and right are identified and depict the equatorial cross-section of $\Sigma$ in this slice, which gets filled in by $B^3$'s in the future and past. 
    
    \begin{figure}[h!]
        \centering
        \includegraphics[width=0.5\linewidth]{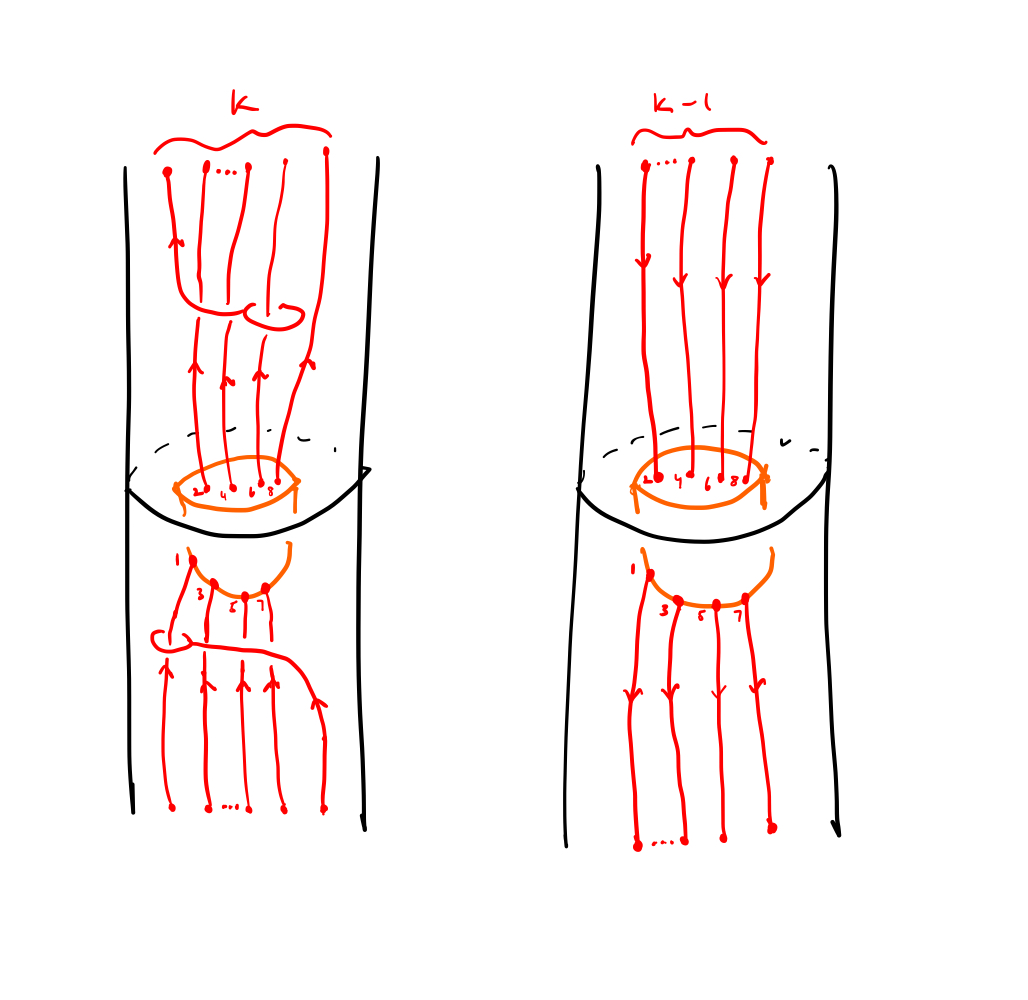}
        \caption{$\beta_k$ in the $t=0$ slice of $L\# R$. The numbering on the orange $S^2$'s denotes which points are identified.}
        \label{beta}
    \end{figure}

    \noindent First we note that $\beta_k$ is not isotopic to the identity. If it were, then after gluing an $S^2 \times D^2$ to $R$ and extending $\beta_k$ over it by the identity, the extension $\hat{\beta_k}$ would also be isotopic to the identity as a diffeomorphism of $L \# (R \cup_{S^1 \times S^2}D^2 \times S^2)=S^1 \times B^3 \# S^4 \approx S^1 \times B^3$. On the other hand, we can see that $\hat{\beta_k} \approx  \delta_k$ of \cite{BG} (\Cref{isotope beta to delta}), which they showed to be nontrivial in $\pi_0(\text{Diff}_\partial(S^1 \times B^3)/\text{Diff}_\partial(B^4))$ for $k\geq 4$.
    
    \begin{figure}[h!]
        \centering
        \includegraphics[width=0.75\linewidth]{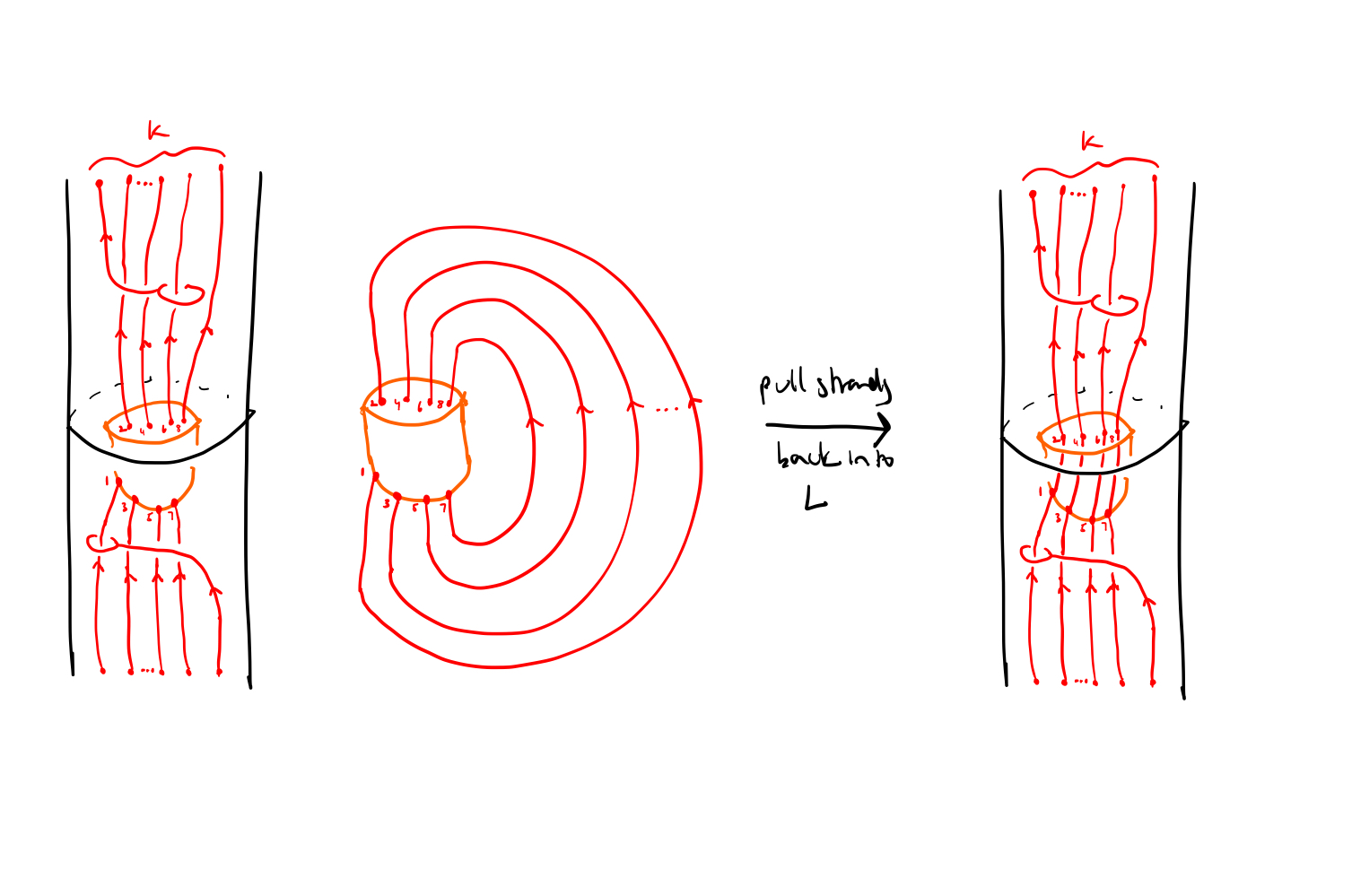}
        \caption{The isotopy of $\hat{\beta_k}$ to $\delta_k.$ Note that we also recover $\delta_k$ if we cap off $L$ rather than $R$.}
        \label{isotope beta to delta}
    \end{figure}

     \noindent Now suppose that $\beta_k(\Sigma)$ is isotopic to $\Sigma$. By \cite{Ce}, $\beta_k$ is isotopic rel $\partial$ to a diffeomorphism $\Bar{\beta}_k$ such that $\Bar{\beta}_k\restriction_\Sigma=id$. Then $\Bar{\beta_k}$ factors as a diffeomorphism $\Bar{\beta_k}=f_L \# f_R$, where $f_L$ and $f_R$ are diffeomorphisms of $L \setminus B^4$ and $R \setminus B^4$ that are the identity near the $S^1 \times S^2 \sqcup S^3$ boundary (\Cref{fLfR}). Moreover, since $\Bar{\beta_k}\not\approx id$ rel $B^4$, then either $f_L\not\approx id$ rel $B^4$ or $f_R\not\approx id$ rel $B^4$. Suppose $f_L\not\approx id$ rel $B^4$. Fill in the missing $B^4$ to $L\setminus B^4$ and extend $f_L$ over it by the identity, and let $\hat{f_L}: S^1 \times B^3 \rightarrow S^1 \times B^3$ denote the extension. Then $\hat{f_L}$ is nontrivial in $\pi_0(\text{Diff}_\partial (S^1 \times B^3)/\text{Diff}_\partial (B^4))$: after gluing an $S^2 \times D^2$ to $R$ in $L \# R$ and extending $f_L \# f_R$ over it by the identity, the extension $\hat{f_L \# f_R}$ agrees with $\hat{f_L}$ outside of a $B^4$ containing $R\setminus B^4 \cup_{S^1\times S^2} D^2 \times S^2 \subset S^1 \times B^3$, so $\hat{f_L}=\hat{f_L \# f_R}$ in $\text{Diff}_\partial (S^1 \times B^3)/\text{Diff}_\partial (B^4)$; since $f_L \# f_R \approx \beta$ in $L \# R$, then $\hat{f_L \# f_R} \approx \hat{\beta}$ in $S^1 \times B^3$, and $\hat{\beta}$ was nontrivial in $\pi_0(\text{Diff}_\partial (S^1 \times B^3)/\text{Diff}_\partial (B^4))$ as we saw above.
     
     \begin{figure}[h!]
         \centering
         \includegraphics[width=0.4\linewidth]{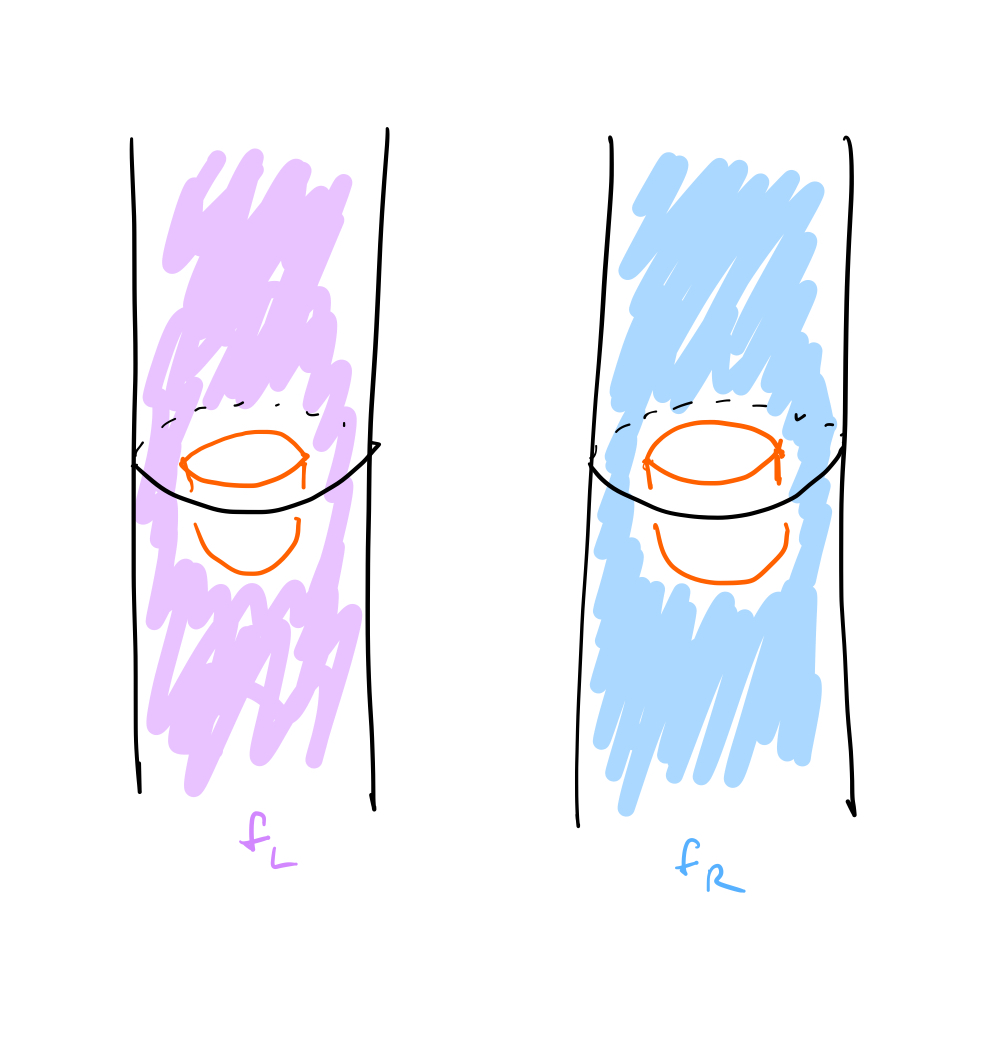}
         \caption{Schematic of $f_L \# f_R$.}
         \label{fLfR}
     \end{figure}

     \noindent Consider the $k$-fold cyclic cover $X\approx \#_{k+1} S^1 \times B^3$ of $L \# R$ obtained by "unwinding" the $R$ factor; see \Cref{2D pic} for the analogous 2-dimensional picture.

     \begin{figure}[h!]
         \centering
         \includegraphics[width=0.4\linewidth]{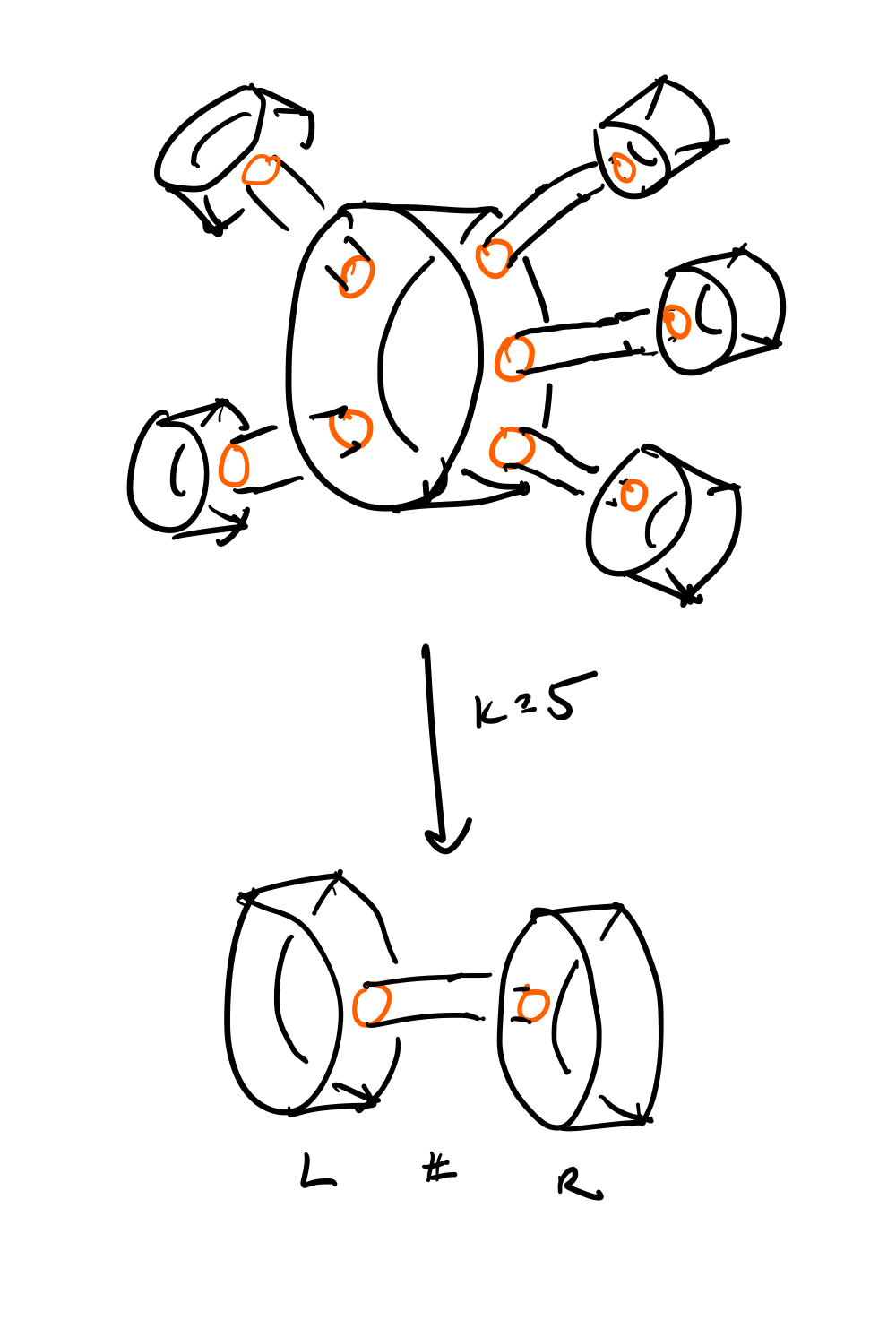}
         \caption{The 2-dimensional analogue of $X$ for $k=5$. The central $S^1 \times B^1$ covers $R$, and the orbiting $S^1 \times B^1$'s cover $L$.}
         \label{2D pic}
     \end{figure}

     \noindent Lifting the barbell representing $\beta_k$ to $X$ and considering the resultant lift $\Tilde{\beta_k}:X \rightarrow X$ of $\beta_k$, we see that $\Tilde{\beta_k}$ is represented by $k$ disjoint copies of a barbell which can be isotoped into a $B^4\subset X$, so that $\Tilde{\beta_k}$ is a composition of trivial elements in $\pi_0(\text{Diff}_\partial(X)/\text{Diff}_\partial(B^4))$ and is therefore trivial too. (\Cref{beta lift})

     \begin{figure}[h!]
         \centering
         \includegraphics[width=0.75\linewidth]{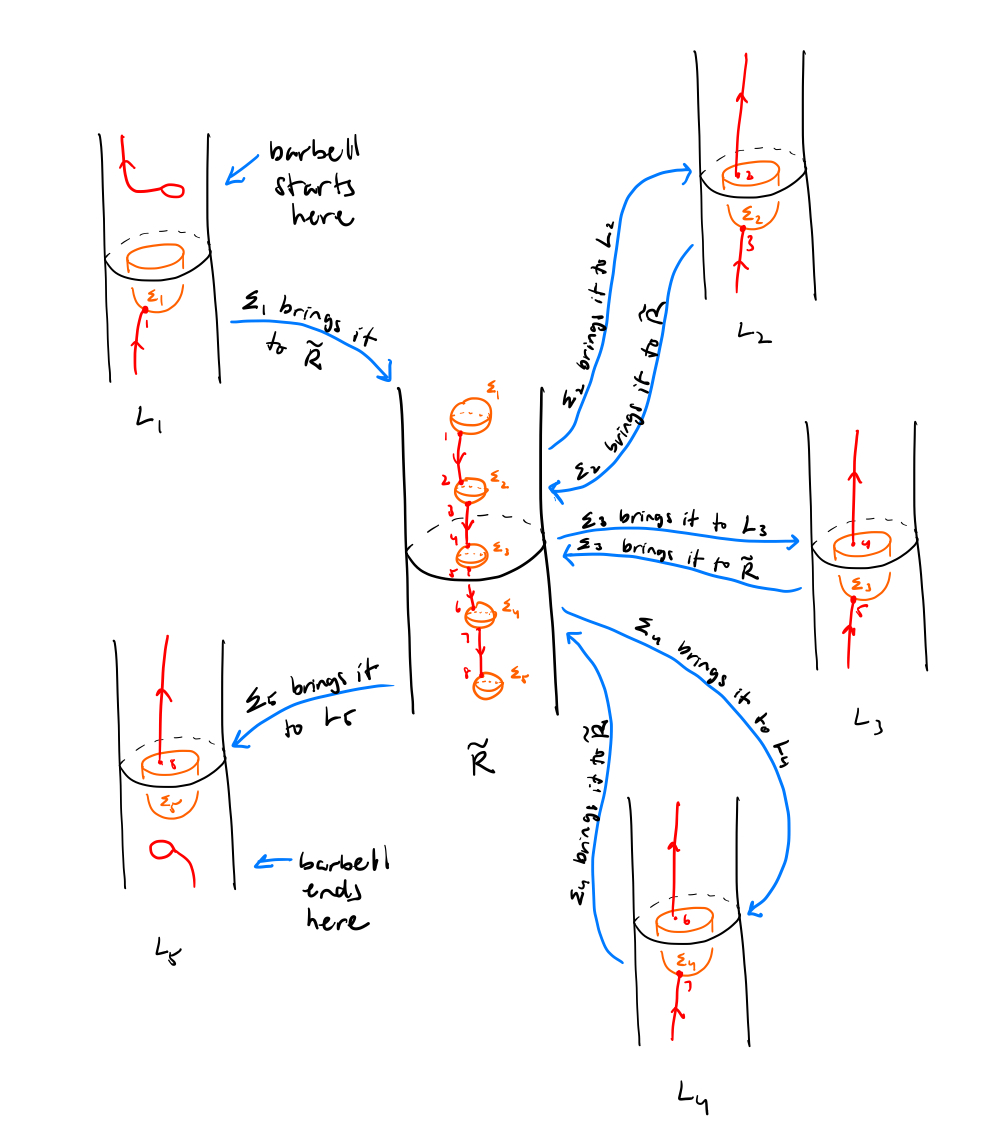}
         \caption{One of the $k$ copies of the barbell for the lift $\Tilde{\beta_k}$ when $k=5$. Here we depict the $t=0$ slice of $X\approx \#_{k+1} S^1 \times D^2 \times [-1,1]$. The orange $S^2$'s are cross-sections of the $k$ $S^3$'s which we identify in the connect sum, and are labeled according to that identification. The barbell begins in $L_1$, hits $\Sigma_1$ and travels to $\Tilde{R}$, then hits $\Sigma_2$ and travels to $L_2$, and so on until ending in $L_5$.  The barbell no longer loops around any $S^1$ factors, so we can pull it into a 4-ball.}
         \label{beta lift}
     \end{figure}

     \noindent On the other hand, consider the lift of $f_L \# f_R$ to $X$; call it $F$. Let $L_1,...,L_k$ denote the $S^1 \times B^3$ summands of $X$ which cover $L$, let $\Tilde{R}$ denote the $S^1 \times B^3$ summand of $X$ which covers $R$, and let $\Sigma_1,...\Sigma_k$ denote the $S^3$'s in $X$ which cover $\Sigma$. Then $F\restriction_{L_i}=f_L$, $F\restriction_{\Tilde{R}}=\Tilde{f_R}$ some lift of $f_R$, and $F$ is the identity on a neighborhood of $\partial X$ and of $\Sigma_1,...,\Sigma_k$.  (\Cref{lift of fL fR})

     \begin{figure}[h!]
         \centering
         \includegraphics[width=0.4\linewidth]{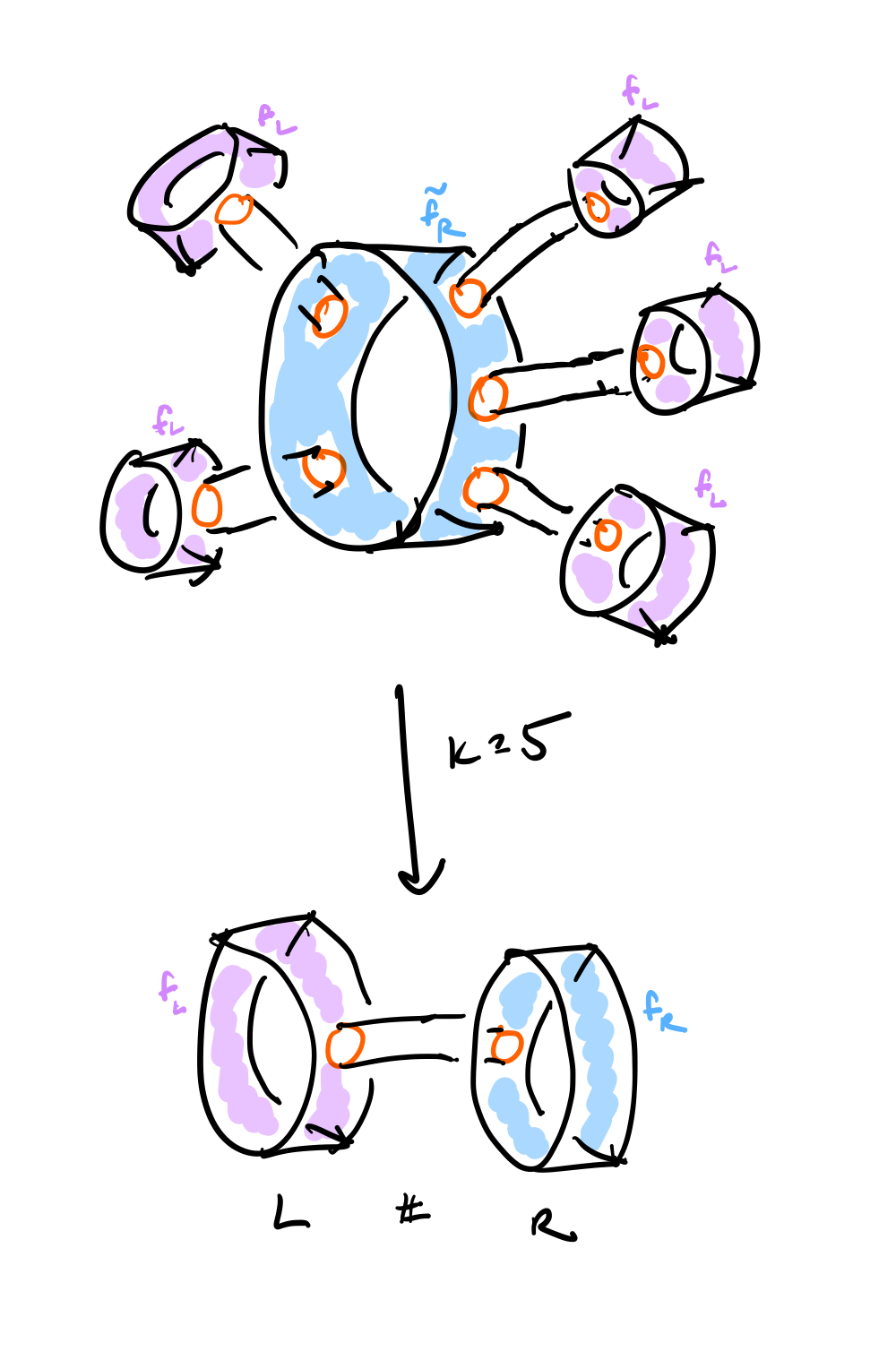}
         \caption{Schematic of the lift $F$ for $k=5$.}
         \label{lift of fL fR}
     \end{figure}

     \noindent But $F$ is nontrivial in $\pi_0(\text{Diff}_\partial(X)/\text{Diff}_\partial(B^4))$; if it were trivial, then after gluing $D^2 \times S^2$'s to $L_2,...,L_k, \Tilde{R}$ and extending $F$ over them by the identity, the extension $\hat{F}:S^1 \times B^3 \rightarrow S^1 \times B^3$ would also be trivial in $\pi_0(\text{Diff}_\partial(S^1 \times B^3)/\text{Diff}_\partial(B^4))$. But $\hat{F}=\hat{f_L}$ outside of a $B^4 \subset S^1 \times B^3$, so $\hat{F}=\hat{f_L}$ in $\text{Diff}_\partial(S^1 \times B^3)/\text{Diff}_\partial(B^4)$, and $\hat{f_L}$ was nontrivial in $\pi_0(\text{Diff}_\partial(S^1 \times B^3)/\text{Diff}_\partial(B^4))$ by assumption. So $F$ is nontrivial in $\pi_0(\text{Diff}_\partial(X)/\text{Diff}_\partial(B^4))$. 

     \noindent Therefore $F\not\approx \Tilde{\beta_k}$, so $f_L \# f_R \not\approx \beta_k$, giving us our contradiction. Therefore $\beta_k(\Sigma)$ is not isotopic to $\Sigma$.

     \noindent Replacing $\beta_k$ by $\beta_i \circ \beta_j^{-1}$ for $i\neq j, i, j \geq 4$, and lifting to the $ij$ -fold cyclic cover of $L\#R$, we see that $\beta_i \circ \beta_j^{-1}(\Sigma)\not\approx \Sigma$; therefore all the $\beta_k$'s for $k\geq 4$ give pairwise non-isotopic $S^3$'s in $S^1 \times B^3 \# S^1 \times B^3$, providing our infinite family of non-isotopic splitting spheres for two unlinked, unknotted $S^2$'s in $S^4$. 
\end{proof}

\end{document}